\documentclass{article}

\usepackage{arxiv}
\usepackage{amssymb}
\usepackage{amsmath}
\usepackage{amsmath,amsfonts}
\usepackage{algorithmic}
\usepackage{algorithm}
\usepackage{array}
\usepackage[caption=false,font=normalsize,labelfont=sf,textfont=sf]{subfig}
\usepackage{textcomp}
\usepackage{stfloats}
\usepackage{url}
\usepackage{verbatim}
\usepackage{graphicx}
\usepackage[utf8]{inputenc} 
\usepackage[T1]{fontenc}    
\usepackage{hyperref}       
\usepackage{url}            
\usepackage{booktabs}       
\usepackage{amsfonts}       
\usepackage{nicefrac}       
\usepackage{microtype}      
\usepackage{cleveref}       
\usepackage{lipsum}         
\usepackage{graphicx}
\usepackage{doi}

\title{The inverse uncertainty distribution of the solutions to a class of higher-order uncertain differential equations}

\author{ {Qiubao Wang} \\
	Department of Mathematics and Physics\\
	Department of Engineering Mechanics\\
	Shijiazhuang Tiedao University \\
	\And
	{Zeman Wang}\thanks{Corresponding Author: wangqiubao12@sina.com} \\
	Department of Mathematics and Physics\\
	Shijiazhuang Tiedao University \\
	\And
	{Zhong Liu}\\
	Department of Mathematics and Physics\\
	Shijiazhuang Tiedao University \\
	\And
	{Zikun Han}\\
	Department of Engineering Mechanics\\
	Shijiazhuang Tiedao University \\
	\And
	{Xiuying Guo}\\
	Department of Mathematics and Physics\\
	Shijiazhuang Tiedao University \\
}

\renewcommand{\shorttitle}{The inverse uncertainty distribution of the solutions to a class of higher-order UDEs}

\begin{document}
\maketitle

\begin{abstract}
	In this paper, we study the higher-order uncertain differential equations (UDEs) as defined by Kaixi Zhang (https://doi.org/10.1007/s10700-024-09422-0),
	mainly focus on the second-order case. We propose a pivotal condition (monotonicity in some sense, see more details in Section 3), introduce the concept of  $\alpha$-paths of UDEs, and demonstrate its properties. Based on this, we derive the inverse uncertainty distribution of the solution. 
\end{abstract}


\section{Introduction}
Events with known frequencies of occurrence are classified as random, while those with unknown frequencies are termed uncertain\textsuperscript{\cite{ref1}}. With the rapid advancement of science and technology, a multitude of uncertain factors have emerged in real life, rendering the phenomena of uncertainty in the objective world undeniable. Consequently, scholars have begun to incorporate these uncertain factors into the establishment of mathematical models, leading to the research and development of uncertainty theory.

Integrating uncertain factors into differential equations results in the formation of UDEs, a type of differential equation established by Liu in 2007\textsuperscript{\cite{ref1}}, designed to describe the dynamics of uncertain phenomena. Yao and Chen provided an effective formula for calculating the inverse uncertainty distribution of the solutions to UDEs, known as the "Yao-Chen Formula"\textsuperscript{\cite{ref2}}.

The Yao-Chen Formula yields a family of solutions to ordinary differential equations, denoted as  $\alpha$-paths, and it has been indicated\textsuperscript{\cite{ref2}} that these $\alpha$-paths represent the inverse uncertainty distribution of the solutions. Therefore, to determine the inverse uncertainty distribution of the solution, we need to solve a family of ordinary differential equations to obtain the $\alpha$-paths. For first-order scalar UDEs, the inverse distribution of the solutions can be determined using the Yao-Chen formula. Building on this foundation, researches have been conducted in many fields, including finance\textsuperscript{\cite{ref3}}, optimal control\textsuperscript{\cite{ref4}}, population growth\textsuperscript{\cite{ref5}}, pharmacokinetics\textsuperscript{\cite{ref6,ref7}}, epidemiology\textsuperscript{\cite{ref8}}, and heat conduction\textsuperscript{\cite{ref9,ref10}}. However, in many practical contexts, first-order scalar UDEs may not fully capture the complexity of real-world scenarios. Often, more intricate cases emerge, such as higher-order, fractional-order, functional differential equations, etc. Thus, investigating theories concerning those higher-order UDEs is of significant importance. This paper primarily focuses on the inverse uncertainty distribution problem of a class of higher-order UDEs.

To determine the inverse uncertainty distribution of the solutions to UDEs, it is essential to identify the $\alpha$-paths. The formation of $\alpha$-paths must satisfy fundamental conditions, primarily concerning the monotonicity with respect to $\alpha$. Clearly, to derive the inverse uncertainty distribution of the solutions to UDEs from 
$\alpha$-paths, these $\alpha$-paths must exhibit monotonicity with respect to $\alpha$. That is, at any given time $t$, the value of $X_t^\alpha$ on an $\alpha$-path should monotonically increase with respect to $\alpha$, as illustrated in Figure 1(a). It is impermissible to encounter a scenario as depicted in Figure 1(b) and Figure 1(c).  The main focus of this paper is to study under what conditions the $\alpha$-paths of UDEs behaves like the one in Figure 1(a), rather than exhibiting the situations shown in (b) and (c).

\begin{figure}[H]
	\centering
	\includegraphics[width=0.8\linewidth]{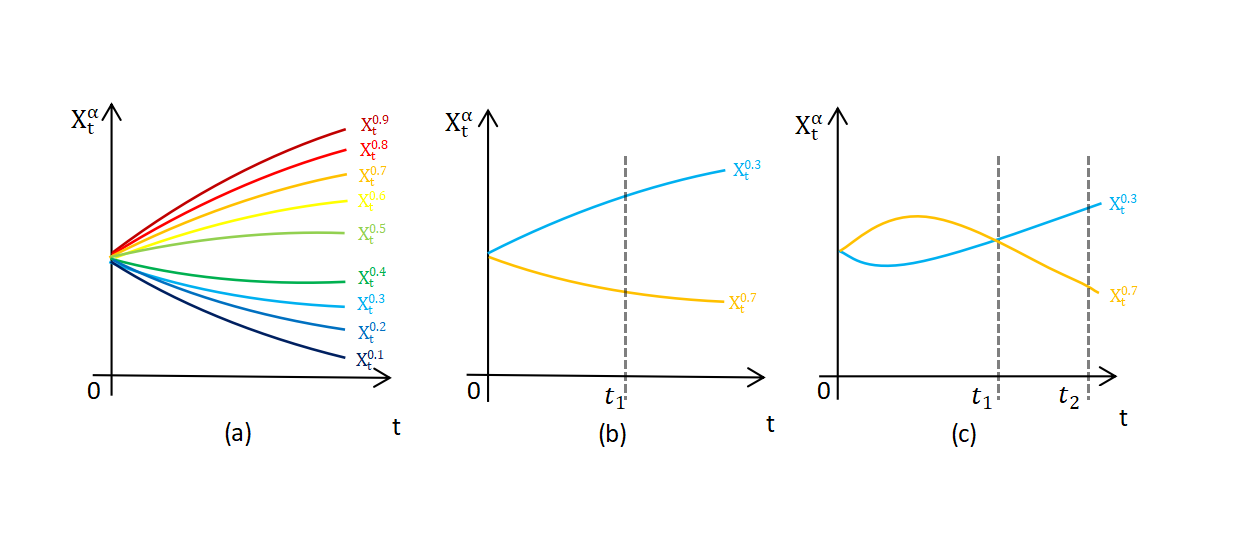}
	\caption{Schematic Diagram}
	\label{fig:1}
\end{figure}

This paper embarks on an exploratory journey for a class of higher-order UDEs, where it becomes imperative to incorporate some restrictive conditions. By ensuring the $\alpha$-paths, we can accurately determine the corresponding inverse distribution of the solutions to UDEs.

The rest of this paper is organized as follows: In Section 2, we reviews related concepts and conclusions of uncertainty theory. In Section 3, we define the $\alpha$-paths for some second-order UDEs and a class of higher-order UDEs, deriving the inverse uncertainty distribution corresponding to the solutions. 
In Section 4, we gives a brief summary to this paper.

\section{Preliminary}

In this section, we introduce some basic concepts and theorems about uncertain processes and uncertain calculus.

\noindent\textbf{Definition 2.1.}(Liu\textsuperscript{\cite{ref11}})  Let $\Gamma$ be a non-empty set, let $\Lambda$ be a $\sigma$-algebra
over $\Gamma$, and let $\mathcal{M}$ be an uncertain measure. Then the triplet $(\Gamma_k,\mathcal{L}_k,\mathcal{M}_k)$ is
called an uncertainty space.

\noindent\textbf{Theorem 2.1.}(Liu\textsuperscript{\cite{ref12}}) (Measure Inversion Theorem) Let $\xi$ be an uncertain
variable with uncertainty distribution $\Phi$. Then for any real number $x$, we have
$$\mathcal{M}\{\xi \le x\} = \Phi(x),
\mathcal{M}\{\xi > x\} = 1-\Phi(x).$$

\noindent\textbf{Definition 2.2.}(Liu\textsuperscript{\cite{ref13}}) An uncertain process $C_t$ is said to be a Liu process if

(i) $C_0 = 0$ and almost all sample paths are Lipschitz continuous,

(ii) $C_t$ has stationary and independent increments,

(iii) every increment $C_{s+t}-C_s$ is a normal uncertain variable with expected value $0$ and variance $t^2$.

\noindent\textbf{Theorem 2.2.}(Yao–Chen Formula\textsuperscript{\cite{ref2}}) Let $X_t$ and $X_t^\alpha$ be the solution and $\alpha$-path of the regular uncertain differential equation
$$	{\rm d}X_t = f(t,X_t){\rm d}t + g(t,X_t){\rm d}C_t,$$
respectively. Then
$$\mathcal{M}\{X_t\le X_t^\alpha,\forall t\} = \alpha,$$
$$\mathcal{M}\{X_t > X_t^\alpha ,\forall t\} = 1-\alpha.$$

\noindent\textbf{Theorem 2.3.}(Zhang\textsuperscript{\cite{ref14}}) Let $f(t,s)$ be an uncertain field, and let $X_t$ and $Y_s$ be general Liu processes. For any partition of the closed region $[0, a] \times [0, b]$ with
$$0 = t_1 < t_2 <\dots < t_{n+1} = a, 0 = s_1 < s_2 < \dots  < s_{m+1} = b,$$
the mesh is written as
\begin{equation*}
	\Delta= \mathop{max}\limits_{
		\substack{1\le i \le n \\  
			1 \le j \le m}}\ \sqrt{(t_{i+1}-t_i)^2+(s_{i+1}-s_i)^2}.
\end{equation*}
Then the double Liu integral of $f(t,s)$ with respect to $X_t$ and $Y_s$ is defined as
$$\iint_{[0,a] \times[0,b]}f(t,s){\rm d}X_t{\rm d}Y_s=\lim_{\Delta \rightarrow 0}\sum_{i = 1}^{n}\sum_{j = 1}^{m}f(t_i,s_i)(X_{t_{i+1}}-X_{t_{i}})(Y_{s_{j+1}}-Y_{s_{j}}),$$
provided that there is an uncertain variable to which the above sum converges almost 
surely as $\Delta \rightarrow 0$. In this case, the uncertain field $f(t,s)$ is said to be integrable with 
respect to $X_t$ and $Y_s$.

\noindent\textbf{Definition 2.3.}(Liu\textsuperscript{\cite{ref13}}; Chen and Ralescu\textsuperscript{\cite{ref15}}; Ye\textsuperscript{\cite{ref16}}) Let $C_t$ be a Liu process, 
and let $Z_t$ be an uncertain process. If there exist two sample-continuous uncertain processes $\mu_t$ and $\delta_t$ such that
$$Z_t=Z_0+\int_{0}^{t} {\mu _s} {\rm d}s+\int_{0}^{t} \delta_s {\rm d}C_s,$$
for any $t\ge 0$, then $Z_t$ is called a general Liu process with drift $\mu_t$ and diffusion $\delta_t$ . 
Furthermore, $Z_t$ has an uncertain differential
$${\rm d}Z_t=\mu_t{\rm d}t+\delta_t {\rm d}C_t,$$
and a first-order derivative
$$\dot{Z}_t=\mu_t+\delta_t\dot{C}_t,$$
where  $\dot{C}_t$ is the formal derivative ${\rm d}C_t/{\rm d}t$.

\noindent\textbf{Theorem 2.4.}(Zhang\textsuperscript{\cite{ref14}})(Fubini's Theorem) Let $X_t$ and $Y_s$ be general Liu processes. Suppose $f(t,s)$
is an integrable uncertain field with respect to $X_t$ and $Y_s$. Then

1.$\int_{0}^{a}f(t,s){\rm d}X_t$ and $\int_{0}^{b}f(t,s){\rm d}Y_s$ exist almost surely;

2.$\int_{0}^{a}f(t,s){\rm d}X_t$ and $\int_{0}^{b}f(t,s){\rm d}Y_s$ are integrable with respect to  $Y_s$ and $X_t$, respectively;

3.$$\iint_{[0,a] \times[0,b]}f(t,s){\rm d}Y_s{\rm d}X_t=\int_{0}^{b}\int_{0}^{a}f(t,s){\rm d}X_t{\rm d}Y_s,a.s.$$

\noindent\textbf{Definition 2.4.}(Zhang\textsuperscript{\cite{ref14}})
Let $C_t$ be a Liu process, and let $Z_t$ be an uncertain process. If there exist sample-continuous uncertain processes $\mu _t$ and $\delta_t$ such that
$$Z_t=Z_0+\dot{Z_0}t+\int_{0}^{t}\int_{0}^{s} {\mu _r} {\rm d}r{\rm d}s+\int_{0}^{t}\int_{0}^{s} \delta_r {\rm d}C_r{\rm d}s$$
for any $t \ge 0$, then $Z_t$ is called a second-order Liu process and has a second-order derivative 
$$\ddot{Z_t}=\mu _t+\delta_t \dot{C_t},$$
where $\ddot{Z_t}$  is the formal second-order derivative ${\rm d}^2 Z_t/{\rm d}t^2.$

\noindent\textbf{Theorem 2.5.}(Liu\textsuperscript{\cite{ref13}}) Let $X_t$
be an uncertain process, and let $C_t$ be a Liu process. 
For any partition of the closed interval $[0, a]$ with
$$0 = t_1 < t_2 < \dots < t_{k+1} = a,$$
the mesh is written as
\begin{equation*}
	\Delta= \mathop{max}\limits_{1 \le i \le k}\ |t_{i+1}-t_i|.
\end{equation*}
Then the Liu integral of $X_t$ with respect to $C_t$ is defined as
$$\int_{0}^{a} X_t {\rm d}C_t = \lim_{\Delta \rightarrow 0}\sum_{i = 1}^{k} X_{t_i}\cdot(C_{t_{i+1} }- C_{t_{i}}),$$
provided that the limit exists almost surely and is finite. In this case, the uncertain	process $X_t$ is said to be integrable.

\noindent\textbf{Definition 2.5.}(Yao-Chen\textsuperscript{\cite{ref2}}) Let $\alpha$ be a number between $0$ and $1$. An uncertain differential equation
$${\rm d}X_t = f(t,X_t){\rm d}t + g(t,X_t){\rm d}C_t$$
is said to have an $\alpha$-path $X_t^\alpha $ if it solves the corresponding ordinary differential equation
$${\rm d}X_t^\alpha= f(t,X^\alpha_t){\rm d}t + |g(t,X^\alpha_t)|\Phi^{-1}(\alpha){\rm d}t $$
where $\Phi^{-1}(\alpha)$ is the inverse standard normal uncertainty distribution, i.e.,
$$\Phi^{-1}(\alpha)=\frac{\sqrt{3}}{\pi}\ln\frac{\alpha}{1-\alpha}.$$

\noindent\textbf{Definition 2.6.}(Liu\textsuperscript{\cite{ref11}}) Suppose $f$ and $g$ are continuous functions. An uncertain differential equation
$${\rm d}X_t = f(t,X_t){\rm d}t + g(t,X_t){\rm d}C_t$$
is said to be regular if
$$g(t,X_t)>0, \forall t>0.$$

\noindent\textbf{Theorem 2.6.}(Liu\textsuperscript{\cite{ref11}})
Let $X_t^\alpha$ be the $\alpha$-path of the regular uncertain differential equation
$${\rm d}X_t = f(t,X_t){\rm d}t + g(t,X_t){\rm d}C_t.$$
Then $X_t^\alpha$ is a continuous and strictly increasing function with respect to $\alpha$ at each time $t>0$ .

\noindent\textbf{Definition 2.7.}(Zhang\textsuperscript{\cite{ref14}})
An uncertain process $X_t$ is called a solution of the higher-order uncertain integral  equation if
\begin{equation*}
	\begin{aligned}
		X_t=&\sum_{i=0}^{n-1}\frac{t^i}{i!}X_0^{i}+\int_{0}^{t}\int_{0}^{t_1}\dots\int_{0}^{t_{n-1}} f(t,X_t,X_t^{(1)},X_t^{(2)},\dots,X_t^{(n-1)}) {\rm d}t_n{\rm d}t_{n-1}\dots{\rm d}t_2{\rm d}t_1\\
		&+\int_{0}^{t}\int_{0}^{t_1}\dots\int_{0}^{t_{n-1}} g(t,X_t,X_t^{(1)},X_t^{(2)},\dots,X_t^{(n-1)}){\rm d}C_{t_n}{\rm d}t_{n-1}\dots{\rm d}t_2{\rm d}t_1,
	\end{aligned}
\end{equation*}
where $C_t$ is a Liu process, and $f$ and $g$ are continuous functions. Equivalently, the  above equation can be simply written as the differential form
$$X_t^{(n)}=f(t,X_t,X_t^{(1)},X_t^{(2)}, \dots,X_t^{(n-1)})+g(t,X_t,X_t^{(1)},X_t^{(2)},\dots,X_t^{(n-1)})\dot{C}_t,$$
that is called a higher-order uncertain diferential equation.

\section{Generalization of Yao-Chen Formula to Some Higher-order UDEs }
Initially, we introduce Lemma 3.1, which serves as a foundation for the proof of Theorems 3.1 to 3.4.

\noindent\textbf{Lemma 3.1.} 
Let $X_t$ satisfies the following higher-order UDE:
$$
\begin{cases}
	X_t^{(n)}&=f(t,X_t,X_t^{(1)},X_t^{(2)}, \dots,X_t^{(n-1)})+g(t,X_t,X_t^{(1)},X_t^{(2)}, \dots,X_t^{(n-1)})\dot{C}_t, \\
	X_t^{(k)}(0)&=X_0^{k},k=1,2,\dots,n-1.
\end{cases}
$$
Then, $X_t$ also satisfies the following integral equation:
\begin{equation*}
	\begin{aligned}
		X_t=&\sum_{k=0}^{n-1}\frac{t^k}{k!}X_0^{k}+\frac{1}{(n-1)!}\int_{0}^{t} (t-s)^{n-1}f(s,X_s^{(1)},X_s^{(2)}, \dots,X_s^{(n-1)}){\rm d}s\\
		&+\frac{1}{(n-1)!}\int_{0}^{t}(t-s)^{n-1}g(s,X_s,X_s^{(1)},X_s^{(2)}, \dots,X_s^{(n-1)}){\rm d}C_s.
	\end{aligned}
\end{equation*}

\noindent\textbf{Proof.} 
Introduce two variables, $Z_t$ and $W_t$, which respectively satisfy the following two equations:
$$ 
\begin{cases}
	Z_t^{(n)}&=0, \\
	Z_t^{(k)}(0)&=Z_0^{k},k=1,2,\dots,n-1.
\end{cases}
$$
$$
\begin{cases}
	W_t^{(n)}&=f(t,W_t^{(1)},W_t^{(2)}, \dots,W_t^{(n-1)})+g(t,W_t,W_t^{(1)},W_t^{(2)}, \dots,W_t^{(n-1)})\dot{C}_t, \\
	W_t^{(k)}(0)&=0.
\end{cases}
$$
It can be easily derived $X_t=Z_t+W_t$.
For $Z_t$, we integrate the both sides directly:
\begin{equation} \label{1}
	\begin{aligned}
		Z_t&=\int_{0}^{t} Z_s^{(1)} {\rm d}s \\
		&=-\int_{0}^{t} Z_s^{(1)} {\rm d}(t-s) \\
		& =-Z_s^{(1)}(t-s)|^t_0-\frac{1}{2}\int_{0}^{t} Z_s^{(2)} {\rm d}(t-s)^2\\
		& =Z_0^{1}t+\frac{1}{2}Z_0^{2}t^2-\frac{1}{3 \times 2}\int_{0}^{t} Z_s^{(3)}(s) {\rm d}(t-s)^3\\
		&\dots \\
		& =Z_0^{1}t+\frac{1}{2}Z_0^{2}t^2+\frac{1}{3 \times 2}Z_0^{3}t^3+\dots +\frac{1}{(n-1)!}Z_0^{n-1}t^{(n-1)}\\
		& =\sum_{k=0}^{n-1}\frac{t^k}{k!}Z_0^{k}.
	\end{aligned}
\end{equation}
For $W_t$, we use intregration by parts repeatedly as follows:
\begin{equation} \label{2}
	\begin{aligned}
		W_t&=\int_{0}^{t} W_s^{(1)} {\rm d} s \\
		&=-\int_{0}^{t} W_s^{(1)} {\rm d}(t-s) \\
		& =-W_s^{(1)}(t-s)|^t_0+\int_{0}^{t} (t-s) W_s^{(2)} {\rm d}s\\
		&=-\frac{1}{2}W_0^{2}t^2+\frac{1}{2}\int_{0}^{t} (t-s)^2 W_s^{(3)} {\rm d}s\\
		&\dots\\
		& =\frac{1}{(n-1)!}\int_{0}^{t} (t-s)^{n-1} f(s,W_s^{(1)},W_s^{(2)}, \dots,W_s^{(n-1)}){\rm d}s\\
		&+\frac{1}{(n-1)!}\int_{0}^{t} (t-s)^{n-1}g(s,W_s,W_s^{(1)},W_s^{(2)}, \dots,W_s^{(n-1)}){\rm d}C_s.
	\end{aligned}
\end{equation}
Thus, adding the eqution (1) and (2), we derive the integral equation satisfied by $X_t$:
\begin{equation*}
	\begin{aligned}
		X_t&=Z_t+ W_t \\
		&=\sum_{k=0}^{n-1}\frac{t^k}{k!}X_0^{k}+\frac{1}{(n-1)!}\int_{0}^{t} (t-s)^{n-1} X_s^{(n)}(s) {\rm d}s\\	&=\sum_{k=0}^{n-1}\frac{t^k}{k!}X_0^{k}+\frac{1}{(n-1)!}\int_{0}^{t} (t-s)^{n-1} f(s,X_s^{(1)},X_s^{(2)}, \dots,X_s^{(n-1)}){\rm d}s\\
		&+\frac{1}{(n-1)!}\int_{0}^{t} (t-s)^{n-1}g(s,X_s,X_s^{(1)},X_s^{(2)}, \dots,X_s^{(n-1)}){\rm d}C_s.
	\end{aligned}
\end{equation*}

\noindent\textbf{Definition 3.1.} Let $\alpha$ be a number between $0$ and $1$. An uncertain differential equation
$$ 
\begin{cases}
	\frac{{\rm d}^2X_t}{{\rm d}t^2} = f(t,X_t,\frac{{\rm d}X_t}{{\rm d}t}) + g(t,X_t,\frac{{\rm d}X_t}{{\rm d}t})\frac{{\rm d}C_t}{{\rm d}t},\\
	X_t\Big| _{t=0}=X_0, \frac{{\rm d}X_t}{{\rm d}t}\Big| _{t=0}=Y_0.
\end{cases}
$$
is said to have an $\alpha$-path $X_t^\alpha$ if it solves the corresponding ordinary differential equation
$$ 
\begin{cases}
	\frac{{\rm d}^2X_t^\alpha}{{\rm d}t^2}= f(t,X^\alpha_t,\frac{{\rm d}X^\alpha_t}{{\rm d}t}) + |g(t,X^\alpha_t,\frac{{\rm d}X_t^\alpha }{{\rm d}t})|\Phi^{-1}(\alpha),\\
	X_t^\alpha\Big| _{t=0}=X_0, \frac{{\rm d}X_t^\alpha}{{\rm d}t}\Big| _{t=0}=Y_0.
\end{cases}
$$
where $\Phi^{-1}(\alpha)$ is the inverse standard normal uncertainty distribution, i.e.,
$$\Phi^{-1}(\alpha)=\frac{\sqrt{3}}{\pi}\ln\frac{\alpha}{1-\alpha}.$$

\noindent\textbf{Definition 3.2.} Suppose $f$ and $g$ are continuous functions. An uncertain differential equation
$$\frac{{\rm d}^{n-1} X_t}{{\rm d}t^{n-1}} = f(t,X_t,\dots,\frac{{\rm d}^{n-1} X_t}{{\rm d}t^{n-1}} ) + g(t,X_t,\dots,\frac{{\rm d}^{n-1} X_t}{{\rm d}t^{n-1}})\frac{{\rm d}C_t}{{\rm d}t}$$
is said to be regular if
$$g(t,X_t,\dots,\frac{{\rm d}^{n-1} X_t}{{\rm d}t^{n-1}})>0, \forall t>0.$$
\noindent\textbf{Theorem 3.1.}
Let $X_t^\alpha$ be the $\alpha$-path of the regular uncertain differential equation
$$\frac{{\rm d}^2X_t}{{\rm d}t^2} = f(t,X_t,\frac{{\rm d}X_t}{{\rm d}t}) + g(t,X_t,\frac{{\rm d}X_t}{{\rm d}t})\frac{{\rm d}C_t}{{\rm d}t}.$$
If the conditions
\begin{equation*}
	\frac{ \partial f(x,y,z)}{\partial y} \ge 0, \frac{ \partial g(x,y,z)}{\partial y} \ge 0 \tag{H}
\end{equation*} 
are met, then $X_t^\alpha$ is a continuous and strictly increasing function with respect to $\alpha$ at each time $t>0$.

\noindent\textbf{Proof.}
Let $\Phi^{-1}$ be the inverse standard normal uncertainty distribution,
and let $\alpha$ and $\beta$ be numbers with $0 < \alpha < \beta < 1$. Write the second-order UDE in the form of a system of equations:
\begin{equation}  \label{3}
	\begin{cases}
		{\rm d}X_t^\alpha=Y_t^\alpha{\rm d}t, \\
		{\rm d}Y_t^\alpha =f(t,X_t^\alpha,Y_t^\alpha) + g(t,X_t^\alpha,Y_t^\alpha)\Phi^{-1}(\alpha),\\
		X_0^\alpha=X_0, Y_0^\alpha=Y_0.
	\end{cases}
\end{equation}
\begin{equation}  \label{4}
	\begin{cases}
		{\rm d}X_t^\beta =Y_t^\beta {\rm d}t, \\
		{\rm d}Y_t^\beta =f(t,X_t^\beta,Y_t^\beta) + g(t,X_t^\beta,Y_t^\beta)\Phi^{-1}(\beta),\\
		X_0^\beta=X_0, Y_0^\beta=Y_0.
	\end{cases}
\end{equation}
Define $\mu$ and $\nu$ as:
$$\mu(T,t)=(T-t)[f(t,X_t^{\alpha},Y_t^{\alpha})+g(t,X_t^{\alpha},Y_t^{\alpha})\Phi^{-1}(\alpha)],$$
$$\nu(T,t)=(T-t)[f(t,X_t^\beta,Y_t^\beta)+g(t,X_t^\beta,Y_t^\beta)\Phi^{-1}(\beta)],$$
Since $g(0,X_0,Y_0)>0$, we have 
$$\mu(T,0)<\nu(T,0)$$
By continuity of $\mu$ and $\nu$, there exists a small number $r > 0$ such that
$$\mu(T,t)<\nu(T,t),\forall t\in[0,r].$$
Thus, by lemma 3.1, we have: 
\begin{equation*}
	\begin{aligned}
		X_T^\alpha
		&=X_0+Y_0T+\int_{0}^{T} (T-t)f(t,X_t^\alpha,Y_t^\alpha) {\rm d}t+\int_{0}^{T}(T-t)g(t,X_t^\alpha,Y_t^\alpha)\Phi^{-1}(\alpha){\rm d}t\\
		&<X_0+Y_0T+\int_{0}^{T} (T-t)f(t,X_t^\beta,Y_t^\beta) {\rm d}t+\int_{0}^{T}(T-t)g(t,X_t^\beta,Y_t^\beta)\Phi^{-1}(\beta){\rm d}t\\
		&=X_T^\beta
	\end{aligned}
\end{equation*}
for any time $T \in (0, r]$.

If for any $t > r$, $X_t^\alpha<X_t^\beta$, the theorem holds. We will prove that by contradiction.

Suppose there exists a time $b > r$ at which $X_t^\alpha$ and $X_t^\beta$ first meet, i.e.,

$$X_b^\alpha=X_b^\beta, X_t^\alpha<X_t^\beta, \forall t \in (0,b).$$

The next phase of our proof will be to compare $Y_t^\alpha$ and  $Y_t^\beta$  in a similar way. 
Write
$$\overline{\mu}(t)=f(t,X_t^{\alpha},Y_t^{\alpha})+g(t,X_t^{\alpha},Y_t^{\alpha})\Phi^{-1}(\alpha),$$
$$\overline{\nu}(t)=f(t,X_t^{\beta},Y_t^{\beta})+g(t,X_t^{\beta},Y_t^{\beta})\Phi^{-1}(\beta).$$
Since  $g(0,X_0,Y_0)>0$, we have 
$$\overline{\mu}(0)<\overline{\nu}(0).$$
By continuity of $\overline{\mu}$ and $\overline{\nu}$, there exists a small number $\overline{r} > 0$ such that
$$\overline{\mu}(t)<\overline{\nu}(t),\forall t\in[0,\overline{r}].$$
Thus,
\begin{equation*}
	\begin{aligned}
		Y_t^\alpha
		&=Y_0+\int_{0}^{t}f(s,X_s^\alpha,Y_s^\alpha) {\rm d}s+\int_{0}^{t}g(s,X_s^\alpha,Y_s^\alpha)\Phi^{-1}(\alpha){\rm d}s\\
		&<Y_0+\int_{0}^{t}f(s,X_s^\beta,Y_s^\beta) {\rm d}s+\int_{0}^{t}g(s,X_s^\beta,Y_s^\beta)\Phi^{-1}(\beta){\rm d}s\\
		&=Y_t^\beta.
	\end{aligned}
\end{equation*}
for any time $t \in (0, \overline{r}]$.

We will prove $Y_t^\alpha<Y_t^\beta$ for $\overline{r}< t < b$ by contradiction.

Suppose there exists a time $\overline{r}< \overline{b} < b$ at which  $Y_t^\alpha$ and $Y_t^\beta$ first meet, i.e.,
$$Y_{\overline{b}}^\alpha=Y_{\overline{b}}^\beta,Y_t^\alpha<Y_t^\beta, \forall t \in (0,\overline{b}).$$
Since (H), we have $$f(\overline{b},X_{\overline{b}}^\alpha,Y_{\overline{b}}^\alpha)<f(\overline{b},X_{\overline{b}}^\beta,Y_{\overline{b}}^\beta),0<g(\overline{b},X_{\overline{b}}^\alpha,Y_{\overline{b}}^\alpha)<g(\overline{b},X_{\overline{b}}^\beta,Y_{\overline{b}}^\beta).$$
Then
$$\overline{\mu}(\overline{b})<\overline{\nu}(\overline{b}).$$
By continuity of $\overline{\mu}$ and $\overline{\nu}$, there exists a time $\overline{a} \in (0, \overline{b})$ such that
$$\overline{\mu}(t)<\overline{\nu}(t),t \in [\overline{a} ,\overline{b}].$$
Thus,
\begin{equation*}
	\begin{aligned}
		Y_{\overline{b}}^\alpha
		&=Y_{\overline{a}}^\alpha+\int_{\overline{a}}^{\overline{b}}f(s,X_s^\alpha,Y_s^\alpha) {\rm d}s+\int_{\overline{a}}^{\overline{b}}g(s,X_s^\alpha,Y_s^\alpha)\Phi^{-1}(\alpha){\rm d}s\\
		&<Y_{\overline{a}}^\beta+\int_{\overline{a}}^{\overline{b}}f(s,X_s^\beta,Y_s^\beta) {\rm d}s+\int_{\overline{a}}^{\overline{b}}g(s,X_s^\beta,Y_s^\beta)\Phi^{-1}(\beta){\rm d}s\\
		&=Y_{\overline{b}}^\beta,
	\end{aligned}
\end{equation*}
which is in contradiction with the assumption $	Y_{\overline{b}}^\alpha=Y_{\overline{b}}^\beta.$ Therefore,
$$Y_{t}^\alpha<Y_{t}^\beta, \forall  0<t<b.$$
Integrate equations(3) and(4)
$$X_b^\alpha=X_0 + \int_{0}^{b}Y_t^\alpha {\rm d}t,$$
$$X_b^\beta=X_0 + \int_{0}^{b}Y_t^\alpha {\rm d}t,$$
we have 
$$X_{b}^\alpha<X_{b}^\beta,$$
which is in contradiction with the assumption $X_{b}^\alpha=X_{b}^\beta.$ The theorem is proved.

By Theorem 3.1, we have established that $X_t$ is a continuous and strictly increasing function with respect to $\alpha$. We wonder whether this collection of $\alpha$-paths can determine the inverse uncertainty distribution of the solutions to second-order UDEs. The following Theorem 3.2 provides a definitive answer.

\noindent\textbf{Theorem 3.2.} Let $X_t$ and $X_t^\alpha$ be the solution and $\alpha$-path of a regular uncertain differential equation
$$\frac{{\rm d}^2X_t}{{\rm d}t^2} = f(t,X_t,\frac{{\rm d}X_t}{{\rm d}t}) + g(t,X_t,\frac{{\rm d}X_t}{{\rm d}t})\frac{{\rm d}C_t}{{\rm d}t}. $$
If the conditions (H) are met,
then we have
$$\mathcal{M}\{X_t\le X_t^\alpha,\forall t\} = \alpha,$$
$$\mathcal{M}\{X_t > X_t^\alpha ,\forall t\} = 1-\alpha.$$

\noindent\textbf{Proof.}
Theorem 14.3 in \cite{ref11} constructs an event $\Lambda_1$ with $\mathcal{M}\{\Lambda_1\}=\alpha,$
and shows that for each $\gamma \in \Lambda_1$, there exists a small number $\delta >0$ such that
\begin{equation}  \label{5}
	\begin{aligned}
		\frac{C_s(\gamma)-C_t(\gamma)}{s-t}< \Phi^{-1}(\alpha -\delta)
	\end{aligned}
\end{equation}
for any times $s$ and $t$ with $s > t$, where $\Phi^{-1}$ is the inverse standard normal uncertainty distribution. Write the second-order UDE in the form of a system of equations:
\begin{equation}  \label{6}
	\begin{cases}
		{\rm d}X_t=Y_t{\rm d}t, \\
		{\rm d}Y_t = f(t,X_t,Y_t){\rm d}t + g(t,X_t,Y_t){\rm d}C_t,\\
		X_{t}\Big| _{t=0}=X_0, Y_{t}\Big| _{t=0}=Y_0.
	\end{cases}
\end{equation}
Define $\lambda$, $\mu$ and $\nu$ as:
$$\lambda(T,t,s)=(T-t)[f(t,X_t(\gamma),Y_t(\gamma))+g(t,X_t(\gamma),Y_t(\gamma))\frac{C_s(\gamma)-C_t(\gamma)}{s-t}],$$
$$\mu(T,t)=(T-t)[f(t,X_t(\gamma),Y_t(\gamma))+g(t,X_t(\gamma),Y_t(\gamma))\Phi^{-1}(\alpha-\delta)],$$
$$\nu(T,t)=(T-t)[f(t,X_t^{\alpha},Y_t^{\alpha})+g(t,X_t^{\alpha},Y_t^{\alpha})\Phi^{-1}(\alpha)].$$
Since $g(0,X_0,Y_0)>0$, we have 
$$\mu(T,0)<\nu(T,0).$$
By continuity of $\mu$ and $\nu$, there exists a small number $r > 0$ such that
$$\mu(T,t)<\nu(T,t),\forall t\in[0,r].$$
By inequality(5), for any time $t \in [0, r]$ and any time $s \in (t, \infty)$, 
we have
\begin{equation*}
	\begin{aligned}
		\lambda(T,t,s)
		&=(T-t)[f(t,X_t(\gamma),Y_t(\gamma))+g(t,X_t(\gamma),Y_t(\gamma))\frac{C_s(\gamma)-C_t(\gamma)}{s-t}] \\
		&<(T-t)[f(t,X_t(\gamma),Y_t(\gamma))+g(t,X_t(\gamma),Y_t(\gamma))\Phi^{-1}(\alpha-\delta)]\\
		&=\mu(T,t)<\nu(T,t).
	\end{aligned}
\end{equation*}
By Lemma 3.1, it follows that
\begin{equation*}
	\begin{aligned}
		X_T(\gamma)
		&=X_0+Y_0T+\int_{0}^{T} (T-t)f(t,X_t(\gamma),Y_t(\gamma)) {\rm d}t+\int_{0}^{T}(T-t)g(t,X_t(\gamma),Y_t(\gamma)){\rm d}C_t(\gamma)\\
		&=X_0+Y_0T+\lim_{\Delta \rightarrow 0}\sum_{i = 1}^{k}\lambda(T,t_{i+1},t_i)(t_{i+1}-t_i)\\
		&\le X_0+Y_0T+\lim_{\Delta \rightarrow 0}\sum_{i = 1}^{k}\mu(T,t_i)(t_{i+1}-t_i)\\
		&=X_0+Y_0T+\int_{0}^{T}\mu(T,t){\rm d}t\\
		&<X_0+Y_0T+\int_{0}^{T}\nu(T,t){\rm d}t\\
		&=X_0+Y_0T+\int_{0}^{T} (T-t)f(t,X_t^\alpha,Y_t^\alpha) {\rm d}t+\int_{0}^{T}(T-t)g(t,X_t^\alpha,Y_t^\alpha)\Phi^{-1}(\alpha){\rm d}t\\
		&=X_T^\alpha
	\end{aligned}
\end{equation*}
for any time $T \in (0, r]$.

We will then prove $X_t(\gamma)<X_t^\alpha$ for all $t>r$ by contradiction.

Suppose there exists a time $b > r$ at which $X_t(\gamma)$ and $X_t^\alpha$ first meet, i.e.,
$$X_b(\gamma)=X_b^\alpha,X_t(\gamma)<X_t^\alpha, \forall t \in (0,b).$$

The next phase of our proof will be compare $Y_t^\alpha$ and  $Y_t(\gamma)$ in a similar way. Write
$$\overline{\lambda}(t,s)=f(t,X_t(\gamma),Y_t(\gamma))+g(t,X_t(\gamma),Y_t(\gamma))\frac{C_s(\gamma)-C_t(\gamma)}{s-t},$$
$$\overline{\mu}(t)=f(t,X_t(\gamma),Y_t(\gamma))+g(t,X_t(\gamma),Y_t(\gamma))\Phi^{-1}(\alpha-\delta),$$
$$\overline{\nu}(t)=f(t,X_t^{\alpha},Y_t^{\alpha})+g(t,X_t^{\alpha},Y_t^{\alpha})\Phi^{-1}(\alpha).$$
Since  $g(0,X_0,Y_0)>0$, we have 
$$\overline{\mu}(0)<\overline{\nu}(0).$$
By continuity of $\overline{\mu}$ and $\overline{\nu}$, there exists a small number $\overline{r} > 0$ such that
$$\overline{\mu}(t)<\overline{\nu}(t),\forall t\in[0,\overline{r}].$$
By inequality(5), for any time $t \in [0, \overline{r}]$ and any time $s \in (t, \infty)$,
we have
\begin{equation*}
	\begin{aligned}
		\overline{\lambda}(t,s)
		&=f(t,X_t(\gamma),Y_t(\gamma))+g(t,X_t(\gamma),Y_t(\gamma))\frac{C_s(\gamma)-C_t(\gamma)}{s-t}\\
		&<f(t,X_t(\gamma),Y_t(\gamma))+g(t,X_t(\gamma),Y_t(\gamma))\Phi^{-1}(\alpha-\delta)\\
		&=\overline{\mu}(t)<\overline{\nu}(t).
	\end{aligned}
\end{equation*}
Thus,
\begin{equation*}
	\begin{aligned}
		Y_t(\gamma)
		&=Y_0+\int_{0}^{t} f(s,X_s(\gamma),Y_s(\gamma)) {\rm d}s+\int_{0}^{t}g(s,X_s(\gamma),Y_s(\gamma)) {\rm d}C_s(\gamma)\\
		&=Y_0+\lim_{\Delta \rightarrow 0}\sum_{i = 1}^{k}\overline{\lambda}(t_{i+1},t_i)(t_{i+1}-t_i)\\
		&\le Y_0+\lim_{\Delta \rightarrow 0}\sum_{i = 1}^{k}\overline{\mu}(t_i)(t_{i+1}-t_i)\\
		&=Y_0+\int_{0}^{t} \overline{\mu}(t) {\rm d}t\\
		&<Y_0+\int_{0}^{t} \overline{\nu}(t) {\rm d}t\\
		&=Y_0+\int_{0}^{t}f(s,X_s^\alpha,Y_s^\alpha) {\rm d}s+\int_{0}^{t}g(s,X_s^\alpha,Y_s^\alpha)\Phi^{-1}(\alpha){\rm d}s\\
		&=Y_t^\alpha
	\end{aligned}
\end{equation*}
for any time $t \in (0, \overline{r}]$.

We will prove $Y_t(\gamma)<Y_t^\alpha$ for all $\overline{r}<t <b$ by contradiction.

Suppose there exists a time $\overline{r}<\overline{b} <b$ at which $Y_t(\gamma)$ and $Y_t^\alpha$ first meet, i.e.,
$$Y_{\overline{b}}(\gamma)=Y_{\overline{b}}^\alpha,Y_t(\gamma)<Y_t^\alpha, \forall t \in (0,\overline{b}).$$
Since (H), 
we have $$f(\overline{b},X_{\overline{b}}(\gamma),Y_{\overline{b}}(\gamma))<f(\overline{b},X_{\overline{b}}^\alpha,Y_{\overline{b}}^\alpha),0<g(\overline{b},X_{\overline{b}}(\gamma),Y_{\overline{b}}(\gamma))<g(\overline{b},X_{\overline{b}}^\alpha,Y_{\overline{b}}^\alpha),$$
then
$$\overline{\mu}(\overline{b})<\overline{\nu}(\overline{b}).$$
By continuity of $\overline{\mu}$ and $\overline{\nu}$, there exists a time $\overline{a} \in (0, \overline{b})$ such that
$$\overline{\mu}(t)<\overline{\nu}(t),t \in [\overline{a} ,\overline{b} ]$$
By inequality(5), for any time $t \in [\overline{a} ,\overline{b}]$ and any time $s \in (t, \infty)$,
we have
\begin{equation*}
	\begin{aligned}
		\overline{\lambda}(t,s)
		&=f(t,X_t(\gamma),Y_t(\gamma))+g(t,X_t(\gamma),Y_t(\gamma))\frac{C_s(\gamma)-C_t(\gamma)}{s-t}\\
		&<f(t,X_t(\gamma),Y_t(\gamma))+g(t,X_t(\gamma),Y_t(\gamma))\Phi^{-1}(\alpha-\delta)\\
		&=\overline{\mu}(t)<\overline{\nu}(t).
	\end{aligned}
\end{equation*}
Thus,
\begin{equation*}
	\begin{aligned}
		Y_{\overline{b}}(\gamma)
		&=Y_{\overline{a}}(\gamma)+\int_{\overline{a}}^{\overline{b}} f(s,X_s(\gamma),Y_s(\gamma)) {\rm d}s+\int_{\overline{a}}^{\overline{b}}g(s,X_s(\gamma),Y_s(\gamma)) {\rm d}C_s(\gamma)\\
		&=Y_{\overline{a}}(\gamma)+\lim_{\Delta \rightarrow 0}\sum_{i = 1}^{k}\overline{\lambda}(t_{i+1},t_i)(t_{i+1}-t_i)\\
		&\le Y_{\overline{a}}(\gamma)+\lim_{\Delta \rightarrow 0}\sum_{i = 1}^{k}\overline{\mu}(t_i)(t_{i+1}-t_i)\\
		&=Y_{\overline{a}}(\gamma)+\int_{\overline{a}}^{\overline{b}} \overline{\mu}(t) {\rm d}s\\
		&<Y_{\overline{a}}^\alpha+\int_{\overline{a}}^{\overline{b}} \overline{\nu}(t) {\rm d}s\\
		&=Y_{\overline{a}}^\alpha+\int_{\overline{a}}^{\overline{b}}f(s,X_s^\alpha,Y_s^\alpha) {\rm d}s+\int_{\overline{a}}^{\overline{b}}g(s,X_s^\alpha,Y_s^\alpha)\Phi^{-1}(\alpha){\rm d}s\\
		&=Y_{\overline{b}}^\alpha,
	\end{aligned}
\end{equation*}
which is in contradiction with $	Y_{\overline{b}}(\gamma)=Y_{\overline{b}}^\alpha.$ Therefore,
$$Y_{t}(\gamma)<Y_{t}^\alpha,\forall 0<t<b.$$
Integrate equations(6),
$$X_b(\gamma)=X_0 + \int_{0}^{b}Y_t(\gamma){\rm d}t,$$
$$X_b^\alpha=X_0 + \int_{0}^{b}Y_t^\alpha {\rm d}t,$$
we have 
$$X_{b}(\gamma)<X_{b}^\alpha,$$
which is in contradiction with $X_{b}(\gamma)=X_{b}^\alpha.$  Therefore,
$$X_{t}(\gamma)<X_{t}^\alpha,\forall t>0.$$
Since $\mathcal{M}\{\Lambda_1\}= \alpha,$ we have
\begin{equation}  \label{7}
	\begin{aligned}
		\mathcal{M}\{X_t\le X_t^\alpha,\forall t\} \ge  \alpha.
	\end{aligned}
\end{equation}
Theorem 14.3 in \cite{ref11} constructs an event $\Lambda_2$ with $\mathcal{M}\{\Lambda_2\}=1-\alpha,$
and shows that for each $\gamma \in \Lambda_2$, there exists a small number $\delta >0$ such that
\begin{equation}  \label{8}
	\begin{aligned}
		\frac{C_s(\gamma)-C_t(\gamma)}{s-t}> \Phi^{-1}(\alpha +\delta)
	\end{aligned}
\end{equation}
for any times $s$ and $t$ with $s > t$. Write
$$\lambda(T,t,s)=(T-t)[f(t,X_t(\gamma),Y_t(\gamma))+g(t,X_t(\gamma),Y_t(\gamma))\frac{C_s(\gamma)-C_t(\gamma)}{s-t}],$$
$$\mu(T,t)=(T-t)[f(t,X_t(\gamma),Y_t(\gamma))+g(t,X_t(\gamma),Y_t(\gamma))\Phi^{-1}(\alpha+\delta)],$$
$$\nu(T,t)=(T-t)[f(t,X_t^{\alpha},Y_t^{\alpha})+g(t,X_t^{\alpha},Y_t^{\alpha})\Phi^{-1}(\alpha)].$$
Since $g(0,X_0,Y_0)>0$, we have 
$$\mu(T,0) > \nu(T,0).$$
By continuity of $\mu $ and $\nu$, there exists a small number $r > 0$ such that
$$\mu(T,t) > \nu(T,t),\forall t\in[0,r].$$
By  inequality(8), for any time $t \in [0, r]$ and any time $s \in (t, \infty)$,
we have
\begin{equation*}
	\begin{aligned}
		\lambda(T,t,s)
		&=(T-t)[f(t,X_t(\gamma),Y_t(\gamma))+g(t,X_t(\gamma),Y_t(\gamma))\frac{C_s(\gamma)-C_t(\gamma)}{s-t}] \\
		&>(T-t)[f(t,X_t(\gamma),Y_t(\gamma))+g(t,X_t(\gamma),Y_t(\gamma))\Phi^{-1}(\alpha+\delta)]\\
		&=\mu(T,t)> \nu(T,t).
	\end{aligned}
\end{equation*}
By Lemma 3.1, it follows that
\begin{equation*}
	\begin{aligned}
		X_T(\gamma)
		&=X_0+Y_0T+\int_{0}^{T} (T-t)f(t,X_t(\gamma),Y_t(\gamma)) {\rm d}t+\int_{0}^{T}(T-t)g(t,X_t(\gamma),Y_t(\gamma)){\rm d}C_t(\gamma)\\
		&=X_0+Y_0T+\lim_{\Delta \rightarrow 0}\sum_{i = 1}^{k}\lambda(T,t_{i+1},t_i)(t_{i+1}-t_i)\\
		&\ge  X_0+Y_0T+\lim_{\Delta \rightarrow 0}\sum_{i = 1}^{k}\mu(T,t_i)(t_{i+1}-t_i)\\
		&=X_0+Y_0T+\int_{0}^{T}\mu(T,t){\rm d}t\\
		&>X_0+Y_0T+\int_{0}^{T}\nu(T,t){\rm d}t\\
		&=X_0+Y_0T+\int_{0}^{T} (T-t)f(t,X_t^\alpha,Y_t^\alpha) {\rm d}t+\int_{0}^{T}(T-t)g(t,X_t^\alpha,Y_t^\alpha)\Phi^{-1}(\alpha){\rm d}t\\
		&=X_T^\alpha
	\end{aligned}
\end{equation*}
for any time $T \in (0, r]$.

We will then prove $X_t(\gamma)<X_t^\alpha$ for all $t>r$ by contradiction.

Suppose there exists a time $b > r$ at which $X_t(\gamma)$ and $X_t^\alpha$ first meet, i.e.,
$$X_b(\gamma)=X_b^\alpha,X_t(\gamma)>X_t^\alpha, \forall t \in (0,b).$$
The next phase of our proof will be compare $Y_t^\alpha$ and  $Y_t(\gamma)$ in a similar way. Write
$$\overline{\lambda}(t,s)=f(t,X_t(\gamma),Y_t(\gamma))+g(t,X_t(\gamma),Y_t(\gamma))\frac{C_s(\gamma)-C_t(\gamma)}{s-t},$$
$$\overline{\mu}(t)=f(t,X_t(\gamma),Y_t(\gamma))+g(t,X_t(\gamma),Y_t(\gamma))\Phi^{-1}(\alpha+\delta),$$
$$\overline{\nu}(t)=f(t,X_t^{\alpha},Y_t^{\alpha})+g(t,X_t^{\alpha},Y_t^{\alpha})\Phi^{-1}(\alpha).$$
Since  $g(0,X_0,Y_0)>0$, we have 
$$\overline{\mu}(0) > \overline{\nu}(0).$$
By continuity of $\overline{\mu} $ and $\overline{\nu}$, there exists a small number $\overline{r} > 0$ such that
$$\overline{\mu}(t)> \overline{\nu}(t),\forall t\in[0,\overline{r}].$$
By inequality(8), for any time $t \in [0, \overline{r}]$ and any time $s \in (t, \infty)$,
we have
\begin{equation*}
	\begin{aligned}
		\overline{\lambda}(t,s)
		&=f(t,X_t(\gamma),Y_t(\gamma))+g(t,X_t(\gamma),Y_t(\gamma))\frac{C_s(\gamma)-C_t(\gamma)}{s-t}\\
		&> f(t,X_t(\gamma),Y_t(\gamma))+g(t,X_t(\gamma),Y_t(\gamma))\Phi^{-1}(\alpha+\delta)\\
		&=\overline{\mu}(t)> \overline{\nu}(t).
	\end{aligned}
\end{equation*}
Thus,
\begin{equation*}
	\begin{aligned}
		Y_t(\gamma)
		&=Y_0+\int_{0}^{t} f(s,X_s(\gamma),Y_s(\gamma)) {\rm d}s+\int_{0}^{t}g(s,X_s(\gamma),Y_s(\gamma)) {\rm d}C_s(\gamma)\\
		&=Y_0+\lim_{\Delta \rightarrow 0}\sum_{i = 1}^{k}\overline{\lambda}(t_{i+1},t_i)(t_{i+1}-t_i)\\
		&\ge Y_0+\lim_{\Delta \rightarrow 0}\sum_{i = 1}^{k}\overline{\mu}(t_i)(t_{i+1}-t_i)\\
		&=Y_0+\int_{0}^{t} \overline{\mu}(t) {\rm d}t\\
		&>Y_0+\int_{0}^{t} \overline{\nu}(t) {\rm d}t\\
		&=Y_0+\int_{0}^{t}f(s,X_s^\alpha,Y_s^\alpha) {\rm d}s+\int_{0}^{t}g(s,X_s^\alpha,Y_s^\alpha)\Phi^{-1}(\alpha){\rm d}s\\
		&=Y_t^\alpha
	\end{aligned}
\end{equation*}
for any time $t \in (0, \overline{r}]$.

We will prove $Y_t(\gamma)>Y_t^\alpha$ for all $\overline{r}<t <b$ by contradiction.

Suppose there exists a time $\overline{r}< \overline{b}< b$ at which $Y_t(\gamma)$ and $Y_t^\alpha$ first meet, i.e.,
$$Y_{\overline{b}}(\gamma)=Y_{\overline{b}}^\alpha,Y_t(\gamma)> Y_t^\alpha, \forall t \in (0,b).$$
Since (H),
we have $$f(\overline{b},X_{\overline{b}}(\gamma),Y_{\overline{b}}(\gamma))>f(\overline{b},X_{\overline{b}}^\alpha,Y_{\overline{b}}^\alpha),g(\overline{b},X_{\overline{b}}(\gamma),Y_{\overline{b}}(\gamma)) > g(\overline{b},X_{\overline{b}}^\alpha,Y_{\overline{b}}^\alpha)>0,$$
then
$$\overline{\mu}(\overline{b})> \overline{\nu}(\overline{b}).$$
By continuity of $\overline{\mu} $ and $\overline{\nu}$, there exists a time $\overline{a} \in (0, \overline{b})$ such that
$$\overline{\mu}(t) > \overline{\nu}(t),t \in [\overline{a} ,\overline{b} ].$$
By inequality(8), for any time $t \in [\overline{a} ,\overline{b}]$ and any time $s \in (t, \infty)$,
we have
\begin{equation*}
	\begin{aligned}
		\overline{\lambda}(t,s)
		&=f(t,X_t(\gamma),Y_t(\gamma))+g(t,X_t(\gamma),Y_t(\gamma))\frac{C_s(\gamma)-C_t(\gamma)}{s-t}\\
		&>f(t,X_t(\gamma),Y_t(\gamma))+g(t,X_t(\gamma),Y_t(\gamma))\Phi^{-1}(\alpha+\delta)\\
		&=\overline{\mu}(t)> \overline{\nu}(t).
	\end{aligned}
\end{equation*}
Thus,
\begin{equation*}
	\begin{aligned}
		Y_{\overline{b}}(\gamma)
		&=Y_{\overline{a}}(\gamma)+\int_{\overline{a}}^{\overline{b}} f(s,X_s(\gamma),Y_s(\gamma)) {\rm d}s+\int_{\overline{a}}^{\overline{b}}g(s,X_s(\gamma),Y_s(\gamma)) {\rm d}C_s(\gamma)\\
		&=Y_{\overline{a}}(\gamma)+\lim_{\Delta \rightarrow 0}\sum_{i = 1}^{k}\overline{\lambda}(t_{i+1},t_i)(t_{i+1}-t_i)\\
		&\ge  Y_{\overline{a}}(\gamma)+\lim_{\Delta \rightarrow 0}\sum_{i = 1}^{k}\overline{\mu}(t_i)(t_{i+1}-t_i)\\
		&=Y_{\overline{a}}(\gamma)+\int_{\overline{a}}^{\overline{b}} \overline{\mu}(t) {\rm d}s\\
		&>Y_{\overline{a}}^\alpha+\int_{\overline{a}}^{\overline{b}} \overline{\nu}(t) {\rm d}s\\
		&=Y_{\overline{a}}^\alpha+\int_{\overline{a}}^{\overline{b}}f(s,X_s^\alpha,Y_s^\alpha) {\rm d}s+\int_{\overline{a}}^{\overline{b}}g(s,X_s^\alpha,Y_s^\alpha)\Phi^{-1}(\alpha){\rm d}s\\
		&=Y_{\overline{b}}^\alpha,
	\end{aligned}
\end{equation*}
which is in contradiction with $	Y_{\overline{b}}(\gamma)=Y_{\overline{b}}^\alpha.$ Therefore,
$$Y_{t}(\gamma)>Y_{t}^\alpha,\forall 0<t<b.$$
Integrate equations(6),
$$X_b(\gamma)=X_0 + \int_{0}^{b}Y_t(\gamma){\rm d}t,$$
$$X_b^\alpha=X_0 + \int_{0}^{b}Y_t^\alpha {\rm d}t,$$
we have 
$$X_b(\gamma)> X_b^\alpha,$$
which is in contradiction with
$X_b(\gamma)= X_b^\alpha$. Therefore,
$$X_{t}(\gamma)> X_{t}^\alpha,\forall t>0.$$
Since $\mathcal{M}\{\Lambda_2\}= 1-\alpha,$ we have
\begin{equation} \label{9}
	\begin{aligned}
		\mathcal{M}\{X_t> X_t^\alpha,\forall t\} \ge  1-\alpha.
	\end{aligned}
\end{equation}
It follows from (7),(9) and 
$$\mathcal{M}\{X_t\le X_t^\alpha,\forall t\}+ \mathcal{M}\{X_t > X_t^\alpha,\forall t\} \le 1,$$
that
$$\mathcal{M}\{X_t\le X_t^\alpha,\forall t\} = \alpha,$$
$$\mathcal{M}\{X_t > X_t^\alpha ,\forall t\} = 1-\alpha.$$ 
The proof is now complete.

Determining the inverse uncertainty distribution of solutions to second-order UDEs is more complex than that for first-order UDEs. In our analysis, we employ Condition H, which is fundamentally derived from the comparison theorem. Readers can refer to publications on the comparison theorem for higher-order differential equations\textsuperscript{\cite{ref17}},  partial order relations, and quasi-monotone functions\textsuperscript{\cite{ref18}}.

We now turn to the discussion of higher-order UDEs. In Definition 2.7, we defined higher-order UDEs and presented them in the form of a system of equations(10). The first component of the solution to system (10) is the solution to the equation as defined in Definition 2.7.
\begin{equation}  \label{10}
	\begin{cases}
		{\rm d}X_t= Y_{t_1}{\rm d}t, \\
		{\rm d}Y_{t_1}= Y_{t_2}{\rm d}t,\\
		\dots\\
		{\rm d}Y_{t_{n-1}}= f(t,X_t,Y_{t_1},\dots, Y_{t_{n-1}}){\rm d}t + g(t,X_t,Y_{t_1},\dots, Y_{t_{n-1}}){\rm d}C_t,\\
		X_t\Big| _{t=0}=X_0, Y_{t_1}\Big| _{t=0}=Y_1,\dots,Y_{t_{n-1}}\Big| _{t=0}=Y_{n-1}.
	\end{cases}
\end{equation}

However, since $f$ is an $(n+1)$-ary function, applying the condition of inverse monotonicity becomes more complex and less straightforward to prove. Yet, if our UDE is formulated as in Theorem 3.3(11), we will proceed to provide a proof for the theorem below.

\noindent\textbf{Definition 3.3.} Let $\alpha$ be a number between $0$ and $1$. An uncertain differential equation
\begin{equation*}
	\begin{cases}
		\frac{{\rm d}^nX_t}{{\rm d}t^n} = f(t,X_t,\dots ,\frac{{\rm d}^nX_t}{{\rm d}t^n}) + g(t,X_t,\dots ,\frac{{\rm d}^nX_t}{{\rm d}t^n})\frac{{\rm d}C_t}{{\rm d}t},\\
		X_{t_0}=X_0, \frac{{\rm d}X_t}{{\rm d}t}\Big| _{t=0}=X_1, \dots, \frac{{\rm d}^n X_t}{{\rm d}t^n}\Big| _{t=0}=X_n .
	\end{cases}
\end{equation*}
is said to have an $\alpha$-path $X_t^\alpha$ if it solves the corresponding ordinary differential equation
\begin{equation*}
	\begin{cases}
		\frac{{\rm d}^nX_t}{{\rm d}t^n} = f(t,X_t,\dots ,\frac{{\rm d}^nX_t}{{\rm d}t^n}) + |g(t,X_t,\dots ,\frac{{\rm d}^nX_t}{{\rm d}t^n})|\Phi^{-1}(\alpha),\\
		X_{t_0}=X_0, \frac{{\rm d}X_t}{{\rm d}t}\Big| _{t=0}=X_1, \dots, \frac{{\rm d}^n X_t}{{\rm d}t^n}\Big| _{t=0}=X_n .
	\end{cases}
\end{equation*}
where $\Phi^{-1}(\alpha)$ is the inverse standard normal uncertainty distribution, i.e.,
$$\Phi^{-1}(\alpha)=\frac{\sqrt{3}}{\pi}\ln\frac{\alpha}{1-\alpha}.$$

\noindent\textbf{Theorem 3.3.} 
Let $X_t^\alpha$ be the $\alpha$-path of a regular uncertain differential equation
\begin{equation}  \label{11}
	\begin{aligned}
		\frac{{\rm d}^n X_t}{{\rm d}t^n} = f(t,X_t) + g(t,X_t)\frac{{\rm d}C_t}{{\rm d}t}.
	\end{aligned}
\end{equation}
Then $X_t^\alpha$ is a continuous and strictly increasing function with respect to $\alpha$ at each time $t>0$.

\noindent\textbf{Proof.}
Let $\Phi^{-1}$ be the inverse standard normal uncertainty distribution,
and let $\alpha$ and $\beta$ be numbers with $0 < \alpha < \beta < 1$. Write 
\begin{equation}  \label{12}
	\begin{cases}
		{\rm d}X_t^\alpha= Y_{t_1}^\alpha {\rm d}t, \\
		{\rm d}Y_{t_1}^\alpha= Y_{t_2}^\alpha{\rm d}t,\\
		\dots\\
		{\rm d}Y_{t_{n-1}}^\alpha= f(t,X^\alpha_t){\rm d}t + g(t,X_t^\alpha)\Phi^{-1}(\alpha),\\
		X_{t}^\alpha\Big| _{t=0}=X_0, Y_{t_1}^\alpha\Big| _{t=0}=Y_1,\dots,Y_{t_{n-1}}^\alpha\Big| _{t=0}=Y_{n-1}.
	\end{cases}
\end{equation}

\begin{equation}  \label{13}
	\begin{cases}
		{\rm d}X_t^\beta = Y_{t_1}^\beta {\rm d}t, \\
		{\rm d}Y_{t_1}^\beta = Y_{t_2}^\beta {\rm d}t,\\
		\dots\\
		{\rm d}Y_{t_{n-1}}^\beta = f(t,X^\beta _t ){\rm d}t + g(t,X_t^\beta )\Phi^{-1}(\beta),\\
		X_{t}^\beta\Big| _{t=0}=X_0, Y_{t_1}^\beta \Big| _{t=0}=Y_1,\dots,Y_{t_{n-1}}^\beta \Big| _{t=0}=Y_{n-1}.
	\end{cases}
\end{equation}
$$\mu(T,t)=(T-t)^{n-1}[f(t,X_t^{\alpha})+g(t,X_t^{\alpha})\Phi^{-1}(\alpha)],$$
$$\nu(T,t)=(T-t)^{n-1}[f(t,X_t^\beta)+g(t,X_t^\beta)\Phi^{-1}(\beta)],$$
$$\mu_1(t)=f(t,X_t^{\alpha})+g(t,X_t^{\alpha})\Phi^{-1}(\alpha),$$
$$\nu_1(t)=f(t,X_t^\beta)+g(t,X_t^\beta)\Phi^{-1}(\beta).$$
Since $g(0,X_0)>0$, we have 
$$\mu(T,0)<\nu(T,0)$$
By continuity of $\mu$ and $\nu$, there exists a small number $r > 0$ such that
$$\mu(T,t)<\nu(T,t),\forall t\in[0,r].$$
By Lemma 3.1, it follows that,
\begin{equation*}
	\begin{aligned}
		X_T^\alpha
		&=X_0+\sum_{k=1}^{n-1}\frac{T^k}{k!}Y_k+\frac{1}{(n-1)!}[\int_{0}^{T} (T-t)^{n-1}f(t,X_t^\alpha) {\rm d}t+\int_{0}^{T}(T-t)^{n-1}g(t,X_t^\alpha)\Phi^{-1}(\alpha){\rm d}t]\\
		&<X_0+\sum_{k=1}^{n-1}\frac{T^k}{k!}Y_k+\frac{1}{(n-1)!}[\int_{0}^{T} (T-t)^{n-1}f(t,X_t^\beta) {\rm d}t+\int_{0}^{T}(T-t)^{n-1}g(t,X_t^\beta)\Phi^{-1}(\beta){\rm d}t]\\
		&=X_T^\beta
	\end{aligned}
\end{equation*}
for any time $T \in (0, r]$.

We will prove $X_t^\alpha < X_t^\beta$ for any $t>r$ by contradiction. 

Suppose there exists a time $b > r$ at which $X_t^\alpha$ and $X_t^\beta$ first meet, i.e.,

$$X_b^\alpha=X_b^\beta,X_t^\alpha<X_t^\beta, \forall t \in (0,b).$$
Since $f(b,X_b^\beta)=f(b,X_b^\alpha),g(b,X_b^\beta)=g(b,X_b^\alpha)>0$,
we have 
$$\mu_1(b)<\nu_1(b).$$
By continuity of $\mu_1$ and $\nu_1$, there exists a time $a \in (0, b)$ such that
$$\mu_1(t)<\nu_1(t),t \in [a,b].$$
Then, we have
$$ \frac{(b-t)^{n-1}}{(n-1)!}\mu_1(t)<\frac{(b-t)^{n-1}}{(n-1)!}\nu_1(t), t \in [a,b).$$
Choose $a_1,b_1$ such that $a<a_1<b_1<b$. Thus, we have
$$ \frac{(b-t)^{n-1}}{(n-1)!}\mu_1(t)<\frac{(b-t)^{n-1}}{(n-1)!}\nu_1(t), t \in [a_1,b_1].$$
Note that any continuous function defined on a closed interval can always reach its minimum. Without loss of generality, suppose that $\nu_1(t)-\mu_1(t)$ reaches its minimum on $[a_1,b_1]$ at $t^\ast \in [a_1,b_1]$. Then, we have the following inequality:
\begin{equation*}
	\begin{aligned}
		&\int_{a}^{b} \frac{(b-t)^{n-1}}{(n-1)!}(\nu_1(t)-\mu_1(t)) {\rm d}t\\
		\geqslant&\int_{a_1}^{b_1} \frac{(b-t)^{n-1}}{(n-1)!}(\nu_1(t)-\mu_1(t)) {\rm d}t \\
		\geqslant&\int_{a_1}^{b_1} \frac{(b-b_1)^{n-1}}{(n-1)!}(\nu_1(t^\ast)-\mu_1(t^\ast)) {\rm d}t\\
		=&\frac{(b-b_1)^{n-1}}{(n-1)!}(\nu_1(t^\ast)-\mu_1(t^\ast))(b_1-a_1)\\
		>&0.
	\end{aligned}
\end{equation*}
Thus,
\begin{equation*}
	\begin{aligned}
		X_b^\alpha 
		&=X_0+\sum_{k=1}^{n-1}\frac{(b-a)^k}{k!}Y_k+\frac{1}{(n-1)!}[\int_{a}^{b} (b-t)^{n-1}f(t,X_t^\alpha) {\rm d}t+\int_{a}^{b}(b-t)^{n-1}g(t,X_t^\alpha)\Phi^{-1}(\alpha){\rm d}t]\\
		&<X_0+\sum_{k=1}^{n-1}\frac{(b-a)^k}{k!}Y_k+\frac{1}{(n-1)!}[\int_{a}^{b} (b-t)^{n-1}f(t,X_t^\beta) {\rm d}t+\int_{a}^{b}(b-t)^{n-1}g(t,X_t^\beta)\Phi^{-1}(\beta){\rm d}t]\\
		&=X_b^\beta,
	\end{aligned}
\end{equation*}
which is in contradiction with
$X_b(\gamma)=X_b^\alpha$. The theorem is proved.

By Theorem 3.3, we have established that $X_t$ is a continuous and strictly increasing function with respect to $\alpha$. We wonder whether this collection of $\alpha$-paths can determine the inverse uncertainty distribution of the solutions to higher-order UDEs. The following Theorem 3.4 provides a definitive answer.

\noindent\textbf{Theorem 3.4.}
Let $X_t$ and $X_t^\alpha$ be the solution and $\alpha$-path of the regular uncertain differential equation
$$\frac{{\rm d}^n X_t}{{\rm d}t^n} = f(t,X_t) + g(t,X_t)\frac{{\rm d}C_t}{{\rm d}t},$$
If the conditions (H) are met, then 
$$\mathcal{M}\{X_t\le X_t^\alpha,\forall t\} = \alpha,$$
$$\mathcal{M}\{X_t > X_t^\alpha ,\forall t\} = 1-\alpha.$$

\noindent\textbf{Proof.}
Theorem 14.3 in \cite{ref11} constructs an event $\Lambda_1$ with $\mathcal{M}\{\Lambda_1\}=\alpha,$
and shows that for each $\gamma \in \Lambda_1$, there exists a small number $\delta >0$ such that
\begin{equation}  \label{14}
	\begin{aligned}
		\frac{C_s(\gamma)-C_t(\gamma)}{s-t}< \Phi^{-1}(\alpha -\delta).
	\end{aligned}
\end{equation}
for any times $s$ and $t$ with $s > t$, where $\Phi^{-1}$ be the inverse standard normal uncertainty distribution. Write the higher-order UDE in the form of a system of equations:
\begin{equation}  \label{15}
	\begin{cases}
		{\rm d}X_t= Y_{t_1} {\rm d}t, \\
		{\rm d}Y_{t_1}= Y_{t_2}{\rm d}t,\\
		\dots\\
		{\rm d}Y_{t_{n-1}}= f(t,X_t){\rm d}t + g(t,X_t){\rm d}C_t,\\
		X_t\Big| _{t=0}=X_0, Y_{t_1}\Big| _{t=0}=Y_1,\dots,Y_{t_{n-1}}\Big| _{t=0}=Y_{n-1}.
	\end{cases}
\end{equation}
$$\lambda(T,t,s)=(T-t)^{n-1}[f(t,X_t(\gamma))+g(t,X_t(\gamma))\frac{C_s(\gamma)-C_t(\gamma)}{s-t}],$$
$$\mu(T,t)=(T-t)^{n-1}[f(t,X_t(\gamma))+g(t,X_t(\gamma))\Phi^{-1}(\alpha-\delta)],$$
$$\nu(T,t)=(T-t)^{n-1}[f(t,X_t^{\alpha})+g(t,X_t^{\alpha})\Phi^{-1}(\alpha)],$$
$$\lambda_1(t,s)=f(t,X_t(\gamma))+g(t,X_t(\gamma))\frac{C_s(\gamma)-C_t(\gamma)}{s-t},$$
$$\mu_1(t)=f(t,X_t(\gamma))+g(t,X_t(\gamma))\Phi^{-1}(\alpha-\delta),$$
$$\nu_1(t)=f(t,X_t^{\alpha})+g(t,X_t^{\alpha})\Phi^{-1}(\alpha).$$
Since $g(0,X_0)>0$, we have 
$$\mu(T,0)<\nu(T,0).$$
By continuity of $\mu$ and $\nu$, there exists a small number $r > 0$ such that
$$\mu(T,t)<\nu(T,t),\forall t\in[0,r].$$
By inequality(14), for any time $t \in [0, r]$ and any time $s \in (t, \infty)$,
we have
\begin{equation*}
	\begin{aligned}
		\lambda(T,t,s)
		&=(T-t)^{n-1}[f(t,X_t(\gamma))+g(t,X_t(\gamma))\frac{C_s(\gamma)-C_t(\gamma)}{s-t}] \\
		&<(T-t)^{n-1}[f(t,X_t(\gamma))+g(t,X_t(\gamma))\Phi^{-1}(\alpha-\delta)\\
		&=\mu(T,t)<\nu(T,t).
	\end{aligned}
\end{equation*}
By Lemma 3.1, it follows that
\begin{equation*}
	\begin{aligned}
		X_T(\gamma)
		&=X_0+\sum_{k=1}^{n-1}\frac{T^k}{k!}Y_k+\frac{1}{(n-1)!}[\int_{0}^{T} (T-t)^{n-1}f(t,X_t(\gamma)) {\rm d}t\\
		&+\int_{0}^{T}(T-t)^{n-1}g(t,X_t(\gamma)){\rm d}C_t(\gamma)]\\
		&=X_0+\sum_{k=1}^{n-1}\frac{T^k}{k!}Y_k+\frac{1}{(n-1)!}\lim_{\Delta \rightarrow 0}\sum_{i = 1}^{k}\lambda(T,t_{i+1},t_i)(t_{i+1}-t_i)\\
		&\le X_0+\sum_{k=1}^{n-1}\frac{T^k}{k!}Y_k+\frac{1}{(n-1)!}\lim_{\Delta \rightarrow 0}\sum_{i = 1}^{k}\mu(T,t_i)(t_{i+1}-t_i)\\
		&=X_0+\sum_{k=1}^{n-1}\frac{T^k}{k!}Y_k+\frac{1}{(n-1)!}\int_{0}^{T} \mu(T,t) {\rm d}t\\
		&<X_0+\sum_{k=1}^{n-1}\frac{T^k}{k!}Y_k+\frac{1}{(n-1)!}\int_{0}^{T} \nu(T,t) {\rm d}t\\
		&=X_0+\sum_{k=1}^{n-1}\frac{T^k}{k!}Y_k+\frac{1}{(n-1)!}\int_{0}^{T} (T-t)^{n-1}f(t,X_t^\alpha) {\rm d}t\\
		&+\frac{1}{(n-1)!}\int_{0}^{T}(T-t)^{n-1}g(t,X_t^\alpha)\Phi^{-1}(\alpha){\rm d}t\\
		&=X_T^\alpha
	\end{aligned}
\end{equation*}
for any time $T \in (0, r]$.

We will prove $X_t(\gamma)<X_t^\alpha$ for all $t > r$ by contradiction. Suppose there exists a time $b > r$ at which $X_t(\gamma)$ and $X_t^\alpha$ first meet, i.e.,
$$X_b(\gamma)=X_b^\alpha,X_t(\gamma)<X_t^\alpha, \forall t \in (0,b).$$
Since $f(b,X_b^\alpha)=f(b,X_b(\gamma)),g(b,X_b^\alpha)=g(b,X_b(\gamma))>0$,
we have 
$$\mu_1(b)<\nu_1(b).$$
By continuity of $\mu_1$ and $\nu_1$, there exists a time $a \in (0, b)$ such that
$$\mu_1(t)<\nu_1(t),t \in [a,b].$$
By inequality(14), for any time $t \in [a,b)$ and any time $s \in (t, \infty)$,
we have
\begin{equation*}
	\begin{aligned}
		(b-t)^{n-1}\lambda_1(t,s)
		&=(b-t)^{n-1}[f(t,X_t(\gamma))+g(t,X_t(\gamma))\frac{C_s(\gamma)-C_t(\gamma)}{s-t}] \\
		&<(b-t)^{n-1}[f(t,X_t(\gamma))+g(t,X_t(\gamma))\Phi^{-1}(\alpha-\delta)\\
		&=(b-t)^{n-1}\mu_1(t)<(b-t)^{n-1}\nu_1(t).
	\end{aligned}
\end{equation*}
By similar process in Theorem 3.3, which indicates that the integral of a non-negative function with a positive interval is also positive, we have the following inequality:
\begin{equation*}
	\begin{aligned}
		X_b(\gamma)
		&=X_0+\sum_{k=1}^{n-1}\frac{(b-a)^k}{k!}Y_k+\frac{1}{(n-1)!}\int_{a}^{b} (b-t)^{n-1}f(t,X_t(\gamma)) {\rm d}t\\
		&+\frac{1}{(n-1)!}\int_{a}^{b}(b-t)^{n-1}g(t,X_t(\gamma)){\rm d}C_t(\gamma)\\
		&=X_0+\sum_{k=1}^{n-1}\frac{(b-a)^k}{k!}Y_k+\frac{1}{(n-1)!}\lim_{\Delta \rightarrow 0}\sum_{i = 1}^{k}(b-t_i)^{n-1}\lambda_1(t_{i+1},t_i)(t_{i+1}-t_i)\\
		&\le X_0+\sum_{k=1}^{n-1}\frac{(b-a)^k}{k!}Y_k+\frac{1}{(n-1)!}\lim_{\Delta \rightarrow 0}\sum_{i = 1}^{k}(b-t_i)^{n-1}\mu_1(t_i)(t_{i+1}-t_i)\\
		&=X_0+\sum_{k=1}^{n-1}\frac{(b-a)^k}{k!}Y_k+\frac{1}{(n-1)!}\int_{a}^{b} (b-t)^{n-1}\mu_1(t) {\rm d}t\\
		&<X_0+\sum_{k=1}^{n-1}\frac{(b-a)^k}{k!}Y_k+\frac{1}{(n-1)!}\int_{a}^{b} (b-t)^{n-1}\nu_1(t) {\rm d}t\\
		&=X_0+\sum_{k=1}^{n-1}\frac{(b-a)^k}{k!}Y_k+\frac{1}{(n-1)!}\int_{a}^{b} (b-t)^{n-1}f(t,X_t^\alpha) {\rm d}t\\
		&+\frac{1}{(n-1)!}\int_{a}^{b}(b-t)^{n-1}g(t,X_t^\alpha)\Phi^{-1}(\alpha){\rm d}t\\
		&=X_b^\alpha
	\end{aligned}
\end{equation*}
which is in contradiction with
$X_b(\gamma)=X_b^\alpha$. Therefore,
$$X_{t}(\gamma)<X_{t}^\alpha,\forall t>0.$$
Since $\mathcal{M}\{\Lambda_1\}= \alpha,$ we have
\begin{equation}  \label{16}
	\begin{aligned}
		\mathcal{M}\{X_t\le X_t^\alpha,\forall t\} \ge  \alpha.
	\end{aligned}
\end{equation}
Theorem 14.3 in \cite{ref1} constructs an event $\Lambda_2$ with $\mathcal{M}\{\Lambda_2\}=1-\alpha,$
and shows that for each $\gamma \in \Lambda_2$, there exists a small number $\delta >0$ such that
\begin{equation}  \label{17}
	\begin{aligned}
		\frac{C_s(\gamma)-C_t(\gamma)}{s-t}> \Phi^{-1}(\alpha +\delta)
	\end{aligned}
\end{equation}
for any times $s$ and $t$ with $s > t$, where $\Phi^{-1}$ be the inverse standard normal uncertainty distribution. Write 
$$\lambda(T,t,s)=(T-t)^{n-1}[f(t,X_t(\gamma))+g(t,X_t(\gamma))\frac{C_s(\gamma)-C_t(\gamma)}{s-t}],$$
$$\mu(T,t)=(T-t)^{n-1}[f(t,X_t(\gamma))+g(t,X_t(\gamma))\Phi^{-1}(\alpha + \delta)],$$
$$\nu(T,t)=(T-t)^{n-1}[f(t,X_t^{\alpha})+g(t,X_t^{\alpha})\Phi^{-1}(\alpha)],$$
$$\lambda_1(t,s)=f(t,X_t(\gamma))+g(t,X_t(\gamma))\frac{C_s(\gamma)-C_t(\gamma)}{s-t},$$
$$\mu_1(t)=f(t,X_t(\gamma))+g(t,X_t(\gamma))\Phi^{-1}(\alpha+\delta),$$
$$\nu_1(t)=f(t,X_t^{\alpha})+g(t,X_t^{\alpha})\Phi^{-1}(\alpha).$$
Since $g(0,X_0)>0$, we have 
$$\mu(T,0)>\nu(T,0)$$
By continuity of $\mu$ and $\nu$, there exists a small number $r > 0$ such that
$$\mu(T,t)>\nu(T,t),\forall t\in[0,r].$$
By inequality(17), for any time $t \in [0, r]$ and any time $s \in (t, \infty)$,
we have
\begin{equation*}
	\begin{aligned}
		\lambda(T,t,s)
		&=(T-t)^{n-1}[f(t,X_t(\gamma))+g(t,X_t(\gamma))\frac{C_s(\gamma)-C_t(\gamma)}{s-t}] \\
		&>(T-t)^{n-1}[f(t,X_t(\gamma))+g(t,X_t(\gamma))\Phi^{-1}(\alpha+\delta)\\
		&=\mu(T,t)>\nu(T,t).
	\end{aligned}
\end{equation*}
Thus,
\begin{equation*}
	\begin{aligned}
		X_T(\gamma)
		&=X_0+\sum_{k=1}^{n-1}\frac{T^k}{k!}Y_k+\frac{1}{(n-1)!}\int_{0}^{T} (T-t)^{n-1}f(t,X_t(\gamma)) {\rm d}t\\
		&+\frac{1}{(n-1)!}\int_{0}^{T}(T-t)^{n-1}g(t,X_t(\gamma)){\rm d}C_t(\gamma)\\
		&=X_0+\sum_{k=1}^{n-1}\frac{T^k}{k!}Y_k+\frac{1}{(n-1)!}\lim_{\Delta \rightarrow 0}\sum_{i = 1}^{k}(T-t_i)^{n-1}\lambda(T,t_{i+1},t_i)(t_{i+1}-t_i)\\
		&\ge X_0+\sum_{k=1}^{n-1}\frac{T^k}{k!}Y_k+\frac{1}{(n-1)!}\lim_{\Delta \rightarrow 0}\sum_{i = 1}^{k}(T-t_i)^{n-1}\mu(T,t_i)(t_{i+1}-t_i)\\
		&=X_0+\sum_{k=1}^{n-1}\frac{T^k}{k!}Y_k+\frac{1}{(n-1)!}\int_{0}^{T} \mu(T,t) {\rm d}t\\
		&>X_0+\sum_{k=1}^{n-1}\frac{T^k}{k!}Y_k+\frac{1}{(n-1)!}\int_{0}^{T} \nu(T,t) {\rm d}t\\
		&=X_0+\sum_{k=1}^{n-1}\frac{T^k}{k!}Y_k+\frac{1}{(n-1)!}\int_{0}^{T} (T-t)^{n-1}f(t,X_t^\alpha) {\rm d}t\\
		&+\frac{1}{(n-1)!}\int_{0}^{T}(T-t)^{n-1}g(t,X_t^\alpha)\Phi^{-1}(\alpha){\rm d}t\\
		&=X_T^\alpha
	\end{aligned}
\end{equation*}
for any time $T \in (0, r].$

We will prove $X_t(\gamma)>X_t^\alpha$ for all $t > r$ by contradiction. 

Suppose there exists a time $b > r$ at which $X_t(\gamma)$ and $X_t^\alpha$ first meet, i.e.,
$$X_b(\gamma)=X_b^\alpha,X_t(\gamma)>X_t^\alpha, \forall t \in (0,b).$$
Since $f(b,X_b^\alpha)=f(b,X_b(\gamma)),g(b,X_b^\alpha)=g(b,X_b(\gamma))>0$,
we have 
$$\mu_1(b)>\nu_1(b).$$
By continuity of $\mu_1$ and $\nu_1$, there exists a time $a \in (0, b)$ such that
$$\mu_1(t)>\nu_1(t),t \in [a,b].$$
By inequality(17), for any time $t \in [a,b]$ and any time $s \in (t, \infty)$,
we have
\begin{equation*}
	\begin{aligned}
		(b-t)^{n-1}\lambda_1(t,s)
		&=(b-t)^{n-1}[f(t,X_t(\gamma))+g(t,X_t(\gamma))\frac{C_s(\gamma)-C_t(\gamma)}{s-t}] \\
		&>(b-t)^{n-1}[f(t,X_t(\gamma))+g(t,X_t(\gamma))\Phi^{-1}(\alpha+\delta)\\
		&=(b-t)^{n-1}\mu_1(t)>(b-t)^{n-1}\nu_1(t).
	\end{aligned}
\end{equation*}
Thus,
\begin{equation*}
	\begin{aligned}
		X_b(\gamma)
		&=X_0+\sum_{k=1}^{n-1}\frac{(b-a)^k}{k!}Y_k+\frac{1}{(n-1)!}\int_{a}^{b} (b-t)^{n-1}f(t,X_t(\gamma)) {\rm d}t\\
		&+\frac{1}{(n-1)!}\int_{a}^{b}(b-t)^{n-1}g(t,X_t(\gamma)){\rm d}C_t(\gamma)\\
		&=X_0+\sum_{k=1}^{n-1}\frac{(b-a)^k}{k!}Y_k+\frac{1}{(n-1)!}\lim_{\Delta \rightarrow 0}\sum_{i = 1}^{k}(b-t)^{n-1}\lambda_1(t_{i+1},t_i)(t_{i+1}-t_i)\\
		&\ge X_0+\sum_{k=1}^{n-1}\frac{(b-a)^k}{k!}Y_k+\frac{1}{(n-1)!}\lim_{\Delta \rightarrow 0}\sum_{i = 1}^{k}(b-t)^{n-1}\mu_1(t_i)(t_{i+1}-t_i)\\
		&=X_0+\sum_{k=1}^{n-1}\frac{(b-a)^k}{k!}Y_k+\frac{1}{(n-1)!}\int_{a}^{b} (b-t)^{n-1}\mu_1(t) {\rm d}t\\
		&>X_0+\sum_{k=1}^{n-1}\frac{(b-a)^k}{k!}Y_k+\frac{1}{(n-1)!}\int_{a}^{b} (b-t)^{n-1}\nu_1(t) {\rm d}t\\
		&=X_0+\sum_{k=1}^{n-1}\frac{(b-a)^k}{k!}Y_k+\frac{1}{(n-1)!}\int_{a}^{b}(b-t)^{n-1}f(t,X_t^\alpha) {\rm d}t\\
		&+\frac{1}{(n-1)!}\int_{a}^{b}(b-t)^{n-1}g(t,X_t^\alpha)\Phi^{-1}(\alpha){\rm d}t\\
		&=X_b^\alpha,
	\end{aligned}
\end{equation*}
which is in contradiction with
$X_b(\gamma)=X_b^\alpha$. Therefore,
$$X_{t}(\gamma)> X_{t}^\alpha,\forall t>0.$$
Since $\mathcal{M}\{\Lambda_2\}= 1-\alpha,$ we have
\begin{equation}  \label{18}
	\begin{aligned}
		\mathcal{M}\{X_t> X_t^\alpha,\forall t\} \ge  1-\alpha.
	\end{aligned}
\end{equation}
It follows from (16),(18) and 
$$\mathcal{M}\{X_t\le X_t^\alpha,\forall t\}+ \mathcal{M}\{X_t > X_t^\alpha,\forall t\} \le 1$$
that
$$\mathcal{M}\{X_t\le X_t^\alpha,\forall t\} = \alpha,$$
$$\mathcal{M}\{X_t > X_t^\alpha ,\forall t\} = 1-\alpha.$$
hold.

\section{Conclusions}
The research presented in this paper is of significant importance and has broad applicability in practical applications. If the correct inverse uncertainty distribution of the solutions to UDEs cannot be obtained, it would be impossible to delve into the study of the behavior of UDEs and their solutions. Consequently, the construction of such UDE models would be meaningless.

This paper makes a specific contribution by identifying a sufficient condition for the inverse uncertainty distribution of the solutions to some second-order UDEs and  a class of higher-order UDEs to be determined by their $\alpha$-paths. Under the constraint of this condition, scholars can continue to explore the statistical properties and other applications of the solutions to UDEs.

Based on our research, there are three areas for potential future exploration: 1. It merits investigation whether the sufficient conditions for determining the inverse distribution of solutions to second-order or higher-order UDEs can be expanded and whether there exist necessary and sufficient conditions. 2. While 'Uncertainty Theory' and this paper discuss the transition from ordinary differential equations to UDEs, our research group is currently working on more complex equations, such as uncertain fractional-order differential equations and uncertain functional differential equations with constant delays. We have conceived a potential sufficient condition for the $\alpha$-paths of systems of UDEs to determine the inverse uncertainty distribution of their solutions. Additionally, uncertain partial differential equations also warrant investigation. 3. There is a multitude of differential equations that satisfy our conditions, with numerous examples provided in the final version under submission. On this basis, more in-depth issues can be explored, such as uncertain dynamics, uncertain bifurcation, and uncertain chaos.


\section*{Disclaimer}
This manuscript is a preprint version, which will be enriched with examples in the final version (currently under submission). Should the paper be accepted, the preprint will be updated accordingly. Interested readers are encouraged to stay tuned for further updates. For any inquiries, please feel free to contact wangqiubao12@sina.com.

\end{document}